\documentclass[12pt]{article}
\usepackage[utf8]{inputenc}
\usepackage{mathrsfs}
\usepackage{fullpage}
\usepackage{amssymb}
\usepackage{amsmath}
\usepackage{epsfig}
\usepackage{algorithm}
\usepackage{algorithmic}
\usepackage{latexsym}
\usepackage{amsmath}
\usepackage{amsfonts}
\usepackage{graphicx}
\usepackage{authblk}
\usepackage[normalem]{ulem}
\usepackage{array}
\usepackage{colortbl}
\usepackage{bm}
\usepackage{color}\usepackage[dvipsnames]{xcolor}
\usepackage{tikz}\usetikzlibrary{patterns}
\usepackage{chessboard}
\usepackage{skak}
\usepackage{array,multirow}
\usepackage[
  bookmarks=false, 
  colorlinks,
  citecolor=brown!70!black,
  linkcolor=brown!80!black,
  urlcolor=blue!70!black,
]{hyperref}
\usepackage{relsize,exscale}

\newcommand{\seqnum}[1]{\href{http://oeis.org/#1}{\underline{#1}}}

\newtheorem{defn}{Definition}
\newtheorem{cor}{Corollary}
\newtheorem{thm}{Theorem}

\newcommand{\comm}[1]{}
 	\definecolor{lightlightgray}{rgb}{0.93, 0.93, 0.93}
 		\definecolor{llightgray}{rgb}{0.87, 0.87, 0.87}

\newcolumntype{C}[1]{>{\centering\arraybackslash }b{#1}}

\newcommand\oeis[1]{\href{https://oeis.org/#1}{#1}}

\DeclareMathAlphabet{\mymathbb}{U}{BOONDOX-ds}{m}{n}

\definecolor{bleuclair}{RGB}{186,183,229}
\definecolor{rougeclair}{RGB}{255,230,231}
\definecolor{bleufonce}{RGB}{59,50,114}
\definecolor{rougefonce}{RGB}{127,0,4}

\let\ge\geqslant

\let\geq\geqslant

\title{Partial Motzkin paths with air pockets of the first kind avoiding peaks, valleys or double rises}

\author[1]{Jean-Luc Baril}
\author[2]{Jos{\'e} L. Ram{\'\i}rez}
\affil[1]{\rm LIB, Universit\'e de Bourgogne \protect\\
  B.P. 47 870, 21078 Dijon Cedex France\protect\\
   {\tt E-mail: barjl@u-bourgogne.fr
   }
}

\affil[2]{\rm Departamento de Matem{\'a}ticas, Universidad  Nacional de Colombia\protect\\
 Bogot{\'a}, Colombia\protect\\
   {\tt E-mail: jlramirezr@unal.edu.co
   }
}
\date{\today}

\date{\today}
\begin{document}

\maketitle

\begin{abstract}
 Motzkin paths with air pockets (MAP) of the first kind are defined as a generalization of Dyck paths with air pockets. They are lattice paths in $\mathbb{N}^2$ starting at the origin made of steps $U=(1,1)$, $D_k=(1,-k)$, $k\geq 1$ and $H=(1,0)$, where two down-steps cannot be consecutive. We enumerate MAP and their prefixes  avoiding peaks (resp. valleys, resp. double rise)  according to the length, the type of the last step, and the height of its end-point.  We express our results using Riordan arrays. Finally, we provide constructive bijections between these paths and restricted Dyck and Motzkin paths.
\end{abstract}

\section{Introduction}
In a recent  paper \cite{bakimava}, the authors introduce, study, and enumerate special classes of lattice paths, called {\it Dyck paths with air pockets} (DAP for short). Such paths are  non empty lattice paths in the first quadrant of $\mathbb{Z}^2$ starting at the origin, and consisting of up-steps $U=(1,1)$ and down-steps $D_k=(1,-k)$, $k\geq 1$, where two down-steps cannot be consecutive. These paths can be viewed as ordinary Dyck paths (i.e., paths in $\mathbb{N}^2$ starting at the origin, ending on the $x$-axis and consisting of $U$ and $D_1(=D)$), where each maximal run of down-steps is condensed into one large down-step. As mentioned in \cite{bakimava}, they also correspond to a stack evolution with (partial) reset operations that cannot be consecutive (see for instance \cite{krin}). The authors enumerate these paths with respect to the length, the type (up or down) of the last step and the height of the end-point. Whenever the last point is on the $x$-axis, they prove that the DAP of length $n$ are in one-to-one correspondence with the peak-less Motzkin paths of length $n-1$.  They also investigate the popularity of many patterns in these paths and they give asymptotic approximations. In a second work \cite{bkmv}, the authors make a study for a generalization of these paths by allowing them to go below the $x$-axis. They call these paths Grand Dyck paths with air pockets (GDAP), and they also yield enumerative results for these paths according to the length and several restrictions on the height. In a third paper, Baril and Barry \cite{BaBa} study two generalizations of DAP by allowing some horizontal steps $H=(1,0)$ with some conditions. They call them {\it Motzkin paths with air pockets of the first and second kind}.

In this paper we study Motzkin paths with air pockets of the first kind, which are defined as Motzkin paths (lattice paths in $\mathbb{N}^2$ starting at the origin and made of $U$, $D$, and $H$), where each maximal run of down-steps is condensed into one large down-step. More precisely, we consider lattice paths in $\mathbb{N}^2$ starting at the origin, consisting of steps $U$, $H$, and $D_k$, $k\geq 1$, where {\it two down-steps cannot be consecutive}.
We denote by $\mathcal{D}$ (resp. $\mathcal{M}$, resp. $\mathcal{MP}$) the set of Dyck paths (resp. Motzkin paths, resp. Motzkin paths with air pockets of the first kind).  Moreover, we denote by $\mathcal{PMP}$ the set of partial Motzkin paths with air pockets (PMAP for short). The MAP of the first kind are enumerated by the sequence \seqnum{A114465}. This sequence also counts the Dyck paths  having no ascents of length 2 that start at an odd level.

Throughout the paper,  we will use the following notations. For $k\geq 0$, we consider the generating function $f_k=f_k(z)$ (resp. $g_k=g_k(z)$, resp. $h_k=h_k(z)$), where the coefficient of $z^n$ in the series expansion is the number of  partial Motzkin paths with air pockets of length $n$  ending at height $k$ with an up-step, (resp. with a down-step,   resp. with a horizontal step $H$). 
We introduce the bivariate generating functions
$$F(u,z)=\sum\limits_{k\geq 0} u^kf_k(z), \quad G(u,z)=\sum\limits_{k\geq 0} u^kg_k(z), \mbox{ and } H(u,z)=\sum\limits_{k\geq 0} u^kh_k(z).$$
For short, we also use the notation $F(u), G(u)$, and $H(u)$ for these functions.

A {\em Riordan array} is an infinite lower triangular matrix whose $k$-th column has generating function
$g(z)f(z)^k$ for all $k \ge 0$, for some formal power series  $g(z)$ and $f(z)$, with $g(0) \neq 0$, $f(0)=0$, and $f'(0)\neq 0$. Such a Riordan array is denoted by  $(g(z),f(z))$. We refer to \cite{Riordan, Riordan2} for more details on  Riordan arrays.  Several authors have used Riordan arrays to study lattice paths; see for example \cite{Sprugnoli2,   RAM, RAM2,  Yang, Yang2,  Yang3}.

The outline of this paper is the following. We present enumerative results for partial Motzkin paths with air pockets of the first kind avoiding peaks (resp. avoiding valleys, resp. avoiding double rises), knowing that a {\it peak} is an occurrence $UD_k$ for some $k\geq 1$, a {\it valley} is an occurrence $D_kU$ for some $k\geq 1$, and a {\it double rise} is an occurrence $UU$. For each avoidance, we provide bivariate generating functions that count the PMAP with respect to the length, the type of the last step (up, down or horizontal step) and the height of the end-point. All these results are obtained algebraically by using the famous kernel method for solving several systems of functional equations. We express our results using Riordan arrays and we deduce closed forms for  PMAP of length $n$ ending at height $k$. Finally, we provide constructive bijections between these paths and some restricted Dyck and Motzkin paths.


\section{Partial peak-less Motzkin paths with air pockets}
 
In this section, we study partial Motzkin paths with air pockets of the first kind avoiding occurrences of $UD_i$ for all $i\geq 1$. 
\subsection{Enumerative results}
 Let $P$ be a length $n$ PMAP ending at height $k\geq 0$ and avoiding the occurrences of $UD_i$ for $i\geq 1$. 
 If the last step of $P$ is $U$, then $k\geq 1$ and we have $P=QU$, where $Q$ is a length $(n-1)$ MAP ending at height $k-1$ and avoiding the peaks ($Q$ can be the empty path). So, we obtain the first relation  $f_k=zf_{k-1}+zg_{k-1}+zh_{k-1}$ for $k\geq 1$, anchored with $f_0=1$ by considering the empty path.  If the last step of $P$ is a down-step $D_i$, $i\geq 1$, then we have $P=QD_i$, where $Q$ is a length $(n-1)$ PMAP ending at height $\ell\geq k+1$ with no up- and down-steps at its end, and with no peaks. So, we obtain the second relation  $g_k=z\sum\limits_{\ell\geq k+1} h_\ell$. If the last step of $P$ is a horizontal step $H$,  then we have $P=QH$, where $Q$ is a length $(n-1)$ PMAP ending at height $k$ with no peaks, which implies that $h_k=z(f_k+g_k+h_k)$ for $k\geq 0$. 
 
 Therefore, we have to solve the following system of equations:
 
\begin{equation}\left\{\begin{array}{l}
f_0=1,\mbox{ and } f_k=zf_{k-1}+zg_{k-1}+zh_{k-1}, \quad k\geq 1,\\
g_k=z\sum\limits_{\ell\geq k+1} h_\ell,\quad k\geq 0,\\
h_k=zf_k+zg_k+zh_k,\quad k\geq 0.\\
\end{array}\right.\label{equ1a}
\end{equation}

Multiplying by $u^k$ the recursions in (\ref{equ1a}) and summing over $k$, we have:
\begin{align*}
F(u)&=1+z\sum\limits_{k\geq 1}u^kf_{k-1}+z\sum\limits_{k\geq 1}u^kg_{k-1} +z\sum\limits_{k\geq 1}u^kh_{k-1}\\
&=1+zuF(u)+zuG(u)+zuH(u),\\
G(u)&=z\sum\limits_{k\geq 0}u^k \Bigl(\sum\limits_{\ell\geq k+1} h_\ell     \Bigr)=z\sum\limits_{k\geq 1}h_k(1+u+\cdots + u^{k-1})\\
&=z\sum\limits_{k\geq 1}\frac{u^k-1}{u-1}h_k=\frac{z}{u-1}(H(u)-H(1)),\\
H(u)&=zF(u)+zG(u)+zH(u).
\end{align*}

Notice that we have $F(1)-H(1)=1$ by considering the difference of the first and third equations. Now, setting $h1:=H(1)$ and solving these functional equations, we obtain

$$F(u)={\frac {{\it h1}\,u{z}^{2}+zu+{z}^{2}-u-z+1}{{u}^{2}z+{z}^{2}-u-z+1}},$$
$$G(u)=-{\frac {z \left( {\it h1}\,uz+z{\it h1}-{\it h1}+z \right) }{{u}^{2}z
+{z}^{2}-u-z+1}},\quad 
H(u)={\frac {z \left( z{\it h1}-u+1 \right) }{{u}^{2}z+{z}^{2}-u-z+1}}.$$

In order to compute $h1$, we use the kernel method (see~\cite{ban, pro}) on $H(u)$. We can write the denominator (which is a polynomial in $u$ of degree 2), as $z(u-r)(u-s)$ with

$$r={\frac {1+\sqrt {-4\,{z}^{3}+4\,{z}^{2}-4\,z+1}}{2z}} \mbox{ and }
s={\frac {1-\sqrt {-4\,{z}^{3}+4\,{z}^{2}-4\,z+1}}{2z}}.$$
Plugging $u=s$ (which has a Taylor expansion at $z=0$) in $H(u)z(u-r)(u-s)$, we obtain the equation 
$zh1-s+1 =0,$ which implies that $$h1=\frac{s-1}{z}.$$
Finally,  after simplifying by the factor $(u-s)$ in the numerators and denominators, we obtain 
 $$F(u)=\frac{r}{r-u},\quad G(u)=\frac{s-1}{r-u},\quad\mbox{ and } \quad H(u)=\frac{1}{r-u},$$
which induces that
$$f_k=[u^k] F(u)=\frac{1}{r^k}, \quad g_k=[u^k] G(u)=\frac{s-1}{r^{k+1}},\quad \mbox{ and }\quad
h_k=[u^k] H(u)=\frac{1}{r^{k+1}}.
$$

\begin{thm} The bivariate generating function for the total number of peak-less PMAP with respect to the  length  and the height of the end-point is given by
$$\mathit{Total}(z,u)=\frac{1}{z(r-u)},$$ and we have
$$[u^k] \mathit{Total}(z,u)=\frac{1}{zr^{k+1}}.$$
Finally, setting $t(n,k)=[z^n][u^k]\mathit{Total}(z,u)$, we have for $n\geq 2$ and $k\geq 1$,  $$t(n,k)=t(n,k-1)+t(n-1,k)-t(n-1,k-2)-t(n-2,k),
$$
and setting $t_n:=t(n,0)$, then we have $$t_n=t_{n-1}+\sum\limits_{k=0}^{n-3}t_{k}t_{n-k-3}+\sum\limits_{k=2}^{n-1}\left(t_k-t_{k-1}\right)t_{n-k-1}.$$
\end{thm}
\noindent {\it Proof.} The first three equalities are immediately deduced from the previous results. Now, let us prove the last equality. Any length $n$ peak-less MAP is of the form ($i$) $HP$ where $P$ is a MAP of length $n-1$, or ($ii$) $U Q HD R$, where $Q,R$ are some MAP such that the length of $Q$ lies into $[0,n-3]$, or ($iii$) $P^\sharp Q$,  where $P^\sharp=UP'D_{i}$, $i\geq 2$, and $P'D_{i-1}$ is a MAP  of length lying into $[2,n-1]$. The number of $P'D_{i-1}$ of a given length $k$ is the total number of peak-less MAP of length $k$ minus the number of peak-less MAP of length $k$ and ending with $H$. Taking into account all these cases, we obtain the result.  
\hfill $\Box$

\begin{cor} The generating function that counts all peak-less PMAP  with respect to the length is given by 
$$\mathit{Total}(z,1)=\frac{1}{z(r-1)}.$$
\end{cor}

The first few terms of the series expansion of $\mathit{Total}(z,1)$ are
$$1 + 2z + 4z^2 + 9z^3 + 22z^4 + 56z^5 + 146z^6 + 388z^7 + 1048z^8 + 2869z^9 + O(x^{10}),$$
which corresponds to the sequence \oeis{A152225} in \cite{oeis} counting  Dyck paths of semilength $n+1$ with no peaks of height $0$ (mod 3) and no valleys of height 2 (mod 3); see \cite{Liu}.

\begin{cor} The generating function that counts the peak-less MAP with respect to the length is given by 
$$\mathit{Total}(z,0)=\frac{1}{zr}.$$
\label{cor2}
\end{cor}
The first few terms of the series expansion of $\mathit{Total}(z,0)$ are  
$$1 + z + z^2 + 2z^3 + 5z^4 + 12z^5 + 29z^6 + 73z^7 + 190z^8 + 505z^9+ O(z^{10}),$$ which corresponds to the sequence \oeis{A152171} in \cite{oeis} counting Dyck paths of length $2n$ with no peaks of height $2$ (mod 3) and no valleys of height 1 (mod 3). In Section 2.2, we will exhibit a constructive bijection between these two classes of paths. 

\medskip

Let $\mathcal{T}$ be the infinite matrix $\mathcal{T}:=[t(n,k)]_{n, k\geq 0}$, where $t(n,k)=[z^n][u^k]\mathit{Total}(z,u)$. The first few rows of the matrix $\mathcal{T}$ are
\[\mathcal{T}=
\left(
\begin{array}{ccccccccc}
 1 & 0 & 0 & 0 & 0 & 0 & 0 & 0 & 0 \\
 1 & 1 & 0 & 0 & 0 & 0 & 0 & 0 & 0 \\
 1 & 2 & 1 & 0 & 0 & 0 & 0 & 0 & 0 \\
 2 & 3 & 3 & 1 & 0 & 0 & 0 & 0 & 0 \\
 5 & 6 & 6 & 4 & 1 & 0 & 0 & 0 & 0 \\
 12 & 15 & 13 & 10 & 5 & 1 & 0 & 0 & 0 \\
 29 & 38 & 33 & 24 & 15 & 6 & 1 & 0 & 0 \\
 73 & 96 & 87 & 63 & 40 & 21 & 7 & 1 & 0 \\
 190 & 248 & 229 & 172 & 110 & 62 & 28 & 8 & 1 \\
\end{array}
\right).\]

\begin{figure}[H]
\centering
\includegraphics[scale=0.85]{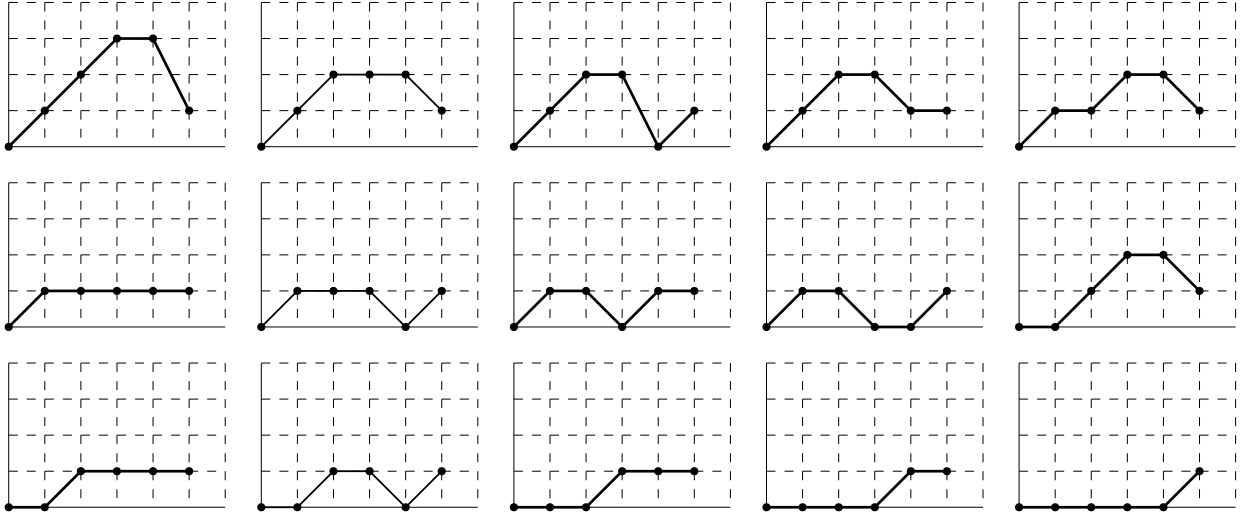}
\caption{Peak-less PMAP of length $5$ ending at height $1$.} \label{fig1}
\end{figure}

In Figure \ref{fig1} we show the peak-less PMAP counted by $t(5,1)=15$. 

\begin{cor}
The matrix $\mathcal{T}=[t(n,k)]_{n, k\geq 0}$ is a Riordan array defined by 
$$\left(C\left(z(1-z+z^2)\right),zC\left(z(1-z+z^2)\right)\right),$$
where $C(z)=\frac{1-\sqrt{1-4z}}{2z}$ is the generating function of the Catalan numbers
$c_n=\frac{1}{n+1}\binom{2n}{n}$.
\end{cor}
\noindent {\it Proof.}
 Indeed, we directly deduce the result from the following. 
\begin{align*}
[u^k]\mathit{Total}(z,u)&=\frac{1}{zr^{k+1}}=\frac{1}{zr}\cdot \frac{1}{r^k}=C\left(z(1-z+z^2)\right)\cdot\left( zC\left(z(1-z+z^2)\right)\right)^k .
\end{align*}
\hfill $\Box$

\begin{cor}\label{corformula}
We have
$$t(n,k)=\sum_{j=0}^{n-k}\frac{k+1}{2(n-j)-k+1}\binom{2(n-j)-k+1}{n-k-j}a(n-k-j,j),$$
where $a(n,k)=(-1)^k\sum_{i=0}^n\binom ni \binom{n - i}{k - 2 i}$.
\end{cor}
\noindent {\it Proof.}
From the Lagrange Inversion Formula (cf. \cite{Lagrange}), we know that 
\begin{align}\label{powCat}
C(z)^k=\sum_{n\geq 0}\frac{k}{2 n + k}\binom{2 n + k}{n}z^n.
\end{align}
From definition of Riordan arrays and Eq. \eqref{powCat} we have
\begin{align*}
t(n,k)&=[z^{n-k}]\left(C\left(z(1-z+z^2)\right)\right)^{k+1}\\
&=[z^{n-k}]\sum_{\ell=0}^{\infty}\frac{k+1}{2\ell+k+1}\binom{2\ell+k+1}{\ell}z^\ell(1-z+z^2)^\ell\\
&=[z^{n-k}]\sum_{\ell=0}^{\infty}\frac{k+1}{2\ell+k+1}\binom{2\ell+k+1}{\ell}z^\ell\sum_{j=0}^{2\ell}a(\ell,j)z^j\\
&=[z^{n-k}]\sum_{j=0}^{\infty}\sum_{\ell=0}^{\infty}\frac{k+1}{2(\ell+\lfloor j/2\rfloor)+k+1}\binom{2(\ell+\lfloor j/2\rfloor)+k+1}{\ell+\lfloor j/2\rfloor}a(\ell+\lfloor j/2\rfloor,j)z^{j+\ell+\lfloor j/2\rfloor}. \end{align*}
If we take $s=j+\ell+\lfloor j/2\rfloor$, then
\begin{align*}
t(n,k)&=[z^{n-k}]\sum_{j=0}^{\infty}\sum_{s=j}^{\infty}\frac{k+1}{2(s-j)+k+1}\binom{2(s-j)+k+1}{s-j}a(s-j,j)z^{s}\\
&=\sum_{j=0}^{n-k}\frac{k+1}{2(n-j)-k+1}\binom{2(n-j)-k+1}{n-k-j}a(n-k-j,j). \end{align*}
\hfill $\Box$

In \cite{Rogers}, Rogers gave an equivalent  characterization of the Riordan arrays. That is, every element not belonging to row $0$ or column $0$ in a Riordan array  can be expressed as a fixed linear combination of the elements in the preceding row. The \emph{$A$-sequence} is defined to be the sequence coefficients of this linear combination. Analogously, Merlini et al. \cite{Merlini}  introduced the \emph{$Z$-sequence}, that characterizes the elements in column $0$, except for the top one. 

An infinite lower triangular matrix $[d_{n,k}]_{n,k\geq0}$ is a Riordan array if and only if  $d_{0,0}\neq 0$ and there exist
 two sequences $(a_0, a_1, a_2, \dots)$, with $a_0\neq 0$, and  $(z_0, z_1, z_2, \dots)$ (called the $A$-sequence and the $Z$-sequence, respectively), such that
\begin{align}\label{AZseq}
d_{n+1,k+1}&=a_0d_{n,k} + a_1d_{n,k+1} + a_2d_{n,k+2} +\cdots & \text{for } n, k\ge0,  \\
d_{n+1,0}&=z_0d_{n,0} + z_1d_{n,1} + z_2d_{n,2} +\cdots & \text{for } n\ge0.
\end{align}
The product of two Riordan arrays $(g(z),f(z))$ and $(h(z),l(z))$ is defined by
\begin{equation}\label{eq:product} (g(z),f(z))*(h(z),l(z))=\left(g(z)h(f(z)), l(f(z))\right). \end{equation}
Under this operation, the set of all Riordan arrays is a group  \cite{Riordan}. The identity element is $I=(1,z)$ and the inverse of $(g(z),f(z))$ is given by
 \begin{equation}\label{invRiordan}
 (g(z), f(z))^{-1}=\left(1/\left(g\circ f^{<-1>}\right)(z), f^{<-1>}(z)\right),
\end{equation}
  where $f^{<-1>}(z)$ denotes the compositional inverse of $f(z)$.

The generating functions for the $A$-sequence and $Z$-sequence of the Riordan array $\mathcal{F}=(g(z),f(z))$,  with inverse $\mathcal{F}^{-1}=(d(z),h(z))$, are given by  (\cite{Merlini, HeS})
\begin{align*}
A(z)=\frac{z}{h(z)} \quad \text{and} \quad Z(z)=\frac{1}{h(z)}\left(1-d_{0,0}d(z)\right),
\end{align*}
respectively.

From the definition of the $A$-sequence and $Z$-sequence for the Riordan arrays we can give an additional recurrence relation for the sequence $t(n,k)$.
\begin{cor}\label{corformula2}
We have
$$t(n+1,k+1)=\sum_{j\geq 0}a(j)t(n,k+j),$$
where $a(n)=(-1)^{n + 1}\sum_{k=1}^n \sum_{j=0}^k \frac{1}{k}\binom{j}{n - k -j}\binom{k}{j}\binom{n - k - 2}{k -1}$ for $n\geq 1$ and $a(0)=1$. Moreover, 
$$t_{n+1}=\sum_{j\geq 0} a(j+1)t(n,j).$$
\end{cor}
\noindent {\it Proof.}
By Equation~\eqref{invRiordan}, the inverse of the matrix $ \mathcal{T}=[t(n,k)]_{n, k\geq 0}$ is given by
$\mathcal{T}^{-1}=\left(g_2(z), zg_2(z)\right)$,
where
$$g_2(z)=\frac{-1 + z^2 + \sqrt{1 - 2 z^2 + 4 z^3 - 3 z^4}}{2z^3}.$$
Therefore, the $A$-sequence and $Z$-sequence of the Riordan array $\mathcal{T}$ have generating functions
\begin{align*}
A(z)=\sum_{n\geq 0}a(n)z^n=\frac{2 z^3}{-1 + z^2 + \sqrt{1 - 2 z^2 + 4 z^3 - 3 z^4}} \quad \text{ and } \quad 
Z(z)=\frac{A(z)-1}{z}.
\end{align*}
The generating function $A(z)$ corresponds with the sequence \seqnum{A247162}, where the explicit formula for $a(n)$ can be found.  From \eqref{AZseq} we obtain the result.
\hfill $\Box$

The first few values of the sequence $a(n)$ for $n\ge0$ are
$$1, \quad 1, \quad 0, \quad 1, \quad 0,\quad  1, \quad -1, \quad 2, \quad -3, \quad 6, \quad -10, \dots.$$

\subsection{A bijective approach}

Corollary~\ref{cor2} proves that the set of peak-less Motzkin paths with air pockets of length $n$ (ending on the $x$-axis) is equinumerous to the set $\mathcal{D}_n(2,1)$ of Dyck paths of length $2n$ with no peak at height 2 (mod 3) and no valley at height 1 (mod 3).

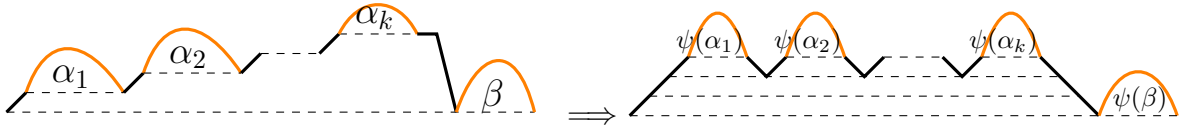
\begin{figure}[H]
\centering
\begin{tikzpicture}[scale=0.26]
\draw[very thick] (20,0) -- (21,1);
\draw[orange, very thick] (21,1) .. controls (22,4) and (24,4) .. (26,1) node[black,xshift=-20,yshift=6] {\large $\alpha_1$};
\draw[dashed] (21,1) -- (26,1);
\draw[dashed] (20,0) -- (47,0);
\draw[very thick] (26,1) -- (27,2);
\draw[orange, very thick] (27,2) .. controls (28,5) and (30,5) .. (32,2) node[black,xshift=-20,yshift=6] {\large $\alpha_2$};
\draw[dashed] (27,2) -- (32,2);
\draw[very thick]  (32,2) -- (33,3);
\draw[dashed] (33,3) -- (36,3);
\draw[very thick]  (36,3) -- (37,4);
\draw[dashed]  (37,4) -- (41,4);
\draw[orange, very thick] (37,4) .. controls (38,6) and (40,6) .. (41,4) node[black,xshift=-16,yshift=6] {\large $\alpha_k$};
\draw[very thick]  (41,4) -- (42,4)--(43,0);
\draw[orange, very thick] (43,0) .. controls (44,3) and (46,4) .. (47,0) node[black,xshift=-16,yshift=6] {\large $\beta$};
\end{tikzpicture}$\quad\Longrightarrow$
\begin{tikzpicture}[scale=0.26]
\draw[orange, very thick] (23,-7) .. controls (24,-4) and (25,-4) .. (26,-7) node[black,xshift=-14,yshift=6] {\smaller $\psi(\alpha_1)$};
\draw[very thick] (20,-10) -- (23,-7);
\draw[very thick] (26,-7) -- (27,-8)-- (28,-7);
 \draw[orange, very thick] (28,-7) .. controls (29,-4) and (30,-4) .. (31,-7) node[black,xshift=-14,yshift=6] {\smaller $\psi(\alpha_2)$};
 \draw[very thick] (31,-7) -- (32,-8)-- (33,-7);
 \draw[dashed] (33,-7) -- (36,-7);
\draw[very thick] (36,-7) -- (37,-8)-- (38,-7);
  \draw[orange, very thick] (38,-7) .. controls (39,-4) and (40,-4) .. (41,-7) node[black,xshift=-14,yshift=6] {\smaller $\psi(\alpha_k)$};
  \draw[very thick] (41,-7) -- (44,-10);
  \draw[orange, very thick] (44,-10) .. controls (45,-7) and (47,-7) .. (48,-10) node[black,xshift=-14,yshift=6] {\smaller$\psi(\beta)$};
 \draw[dashed] (23,-7) -- (26,-7); \draw[dashed] (28,-7) -- (31,-7);\draw[dashed] (38,-7) -- (41,-7);
\draw[dashed] (21,-9) -- (43,-9);\draw[dashed] (22,-8) -- (42,-8);
\draw[dashed] (20,-10) -- (48,-10);
\end{tikzpicture}
\caption{Illustration of the map $\psi$ for the more general case ($iii$) in Definition~\ref{Def1}.}
\label{fig2}
\end{figure}

 Any non-empty peak-less Motzkin path with air pockets is either of the form (1) $H\alpha$  or (2) $U\alpha_1U\alpha_2\cdots U\alpha_k HD_k\beta$,  where $k\geq 1$ and $\alpha,\alpha_1, \ldots, \alpha_k,\beta$ are possibly empty peak-less MAP.  We refer to the left part  of  Figure \ref{fig2} for an illustration of this form.   

\begin{figure}[H]
\centering
\includegraphics[scale=0.75]{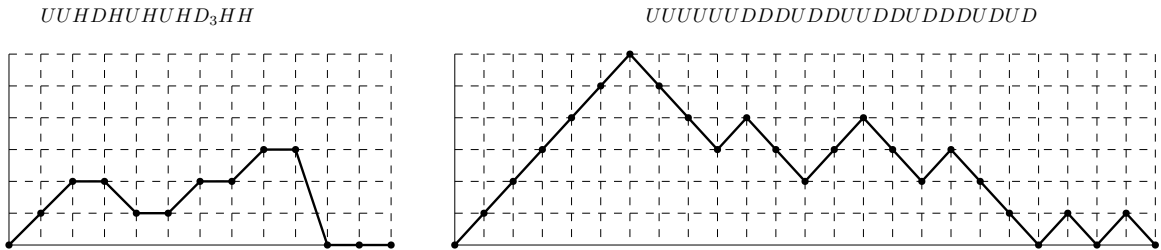}
\caption{$\psi(UUHDHUHUHD_3HH)=UUUUUUDDDUDDUUDDUDDDUDUD$.} \label{fig3}
\end{figure}

Remember that $\mathcal{MP}$ denotes the set of  Motzkin paths with air pockets of the first kind.
\begin{defn} We  recursively define the map $\psi$ from  $\mathcal{MP}$ to $\cup_{n\geq 0}\mathcal{D}_n(2,1)$ as follows. For $\alpha\in\mathcal{MP}$, we set:
$$
\psi(P)=\left\{\begin{array}{llr}
\epsilon&\text{if }P=\epsilon&(i)\\
UD\psi(\alpha)&\text{if }P=H\alpha \text{ with } \alpha\in\mathcal{MP}&(ii)\\
U^3\psi(\alpha_1)DU\psi(\alpha_2)DU\ldots DU\psi(\alpha_k)D^3\psi(\beta)&\text{if }P=U\alpha_1U\alpha_2\ldots U\alpha_k HD_k\beta\text{ with}\\ & k\geq 1 \mbox{ and }\alpha_1,\ldots,\alpha_k,\beta \in\mathcal{MP}. &(iii)
\end{array}\right.
$$
\label{Def1}
\end{defn}
We refer to Figure~\ref{fig3} for an illustration of the third case of the definition of $\psi$.

Due to the recursive definition, the image  of peak-less MAP of length $n$ under $\psi$ is a Dyck path of length $2n$. Moreover it is clear that the obtained path is a Dyck path in $\cup_{n\geq 0}\mathcal{D}_n(2,1)$. For instance (see Figure \ref{fig3} for an illustration of this example).
\begin{align*}
\psi(UUHDHUHUHD_3HH)&=\psi(U\cdot \overbrace{UHDH}^{\alpha_1}\cdot U \cdot \overbrace{H}^{\alpha_2} \cdot U\cdot \overbrace{\epsilon}^{\alpha_3}\cdot H D_3 \cdot \overbrace{HH}^{\beta}\\
&=U^3\psi(UHDH)\cdot DU \cdot \psi(H) \cdot DU \cdot \psi(\epsilon)\cdot  D^3 \cdot \psi(HH)\\
&=U^3\cdot (U^3 D^3 \psi(H)) \cdot DU \cdot UD \cdot DU \cdot  D^3 \cdot UDUD\\
&=U^6D^3UD^2U^2D^2UD^2DUDUD.
\end{align*}

\begin{thm}
For all $n\geq 0$, the map $\psi$ induces a bijection between $\mathcal{MP}_n$ and $\mathcal{D}_{n}(2,1)$.
\label{bij} 
\end{thm}
\noindent {\it Proof.}  Since $\mathcal{MP}_n$ and $\mathcal{D}_n(2,1)$ have the same cardinality (due to Corollary~2 and \oeis{A152171} in \cite{oeis}), it suffices to prove that for $P,Q\in \mathcal{MP}$, $P\neq Q$ implies $\psi(P)\neq \psi(Q)$. A simple induction on $n$ allows to obtain the result.
\hfill $\Box$


\section{Partial valley-less Motzkin paths with air pockets}

Using the same arguments we used for the system of the previous section, we study partial Motzkin paths with air pockets of the first kind avoiding occurrences of $D_iU$ for all $i\geq 1$. 

\subsection{Enumerative results}

In the same way as we done in Section 2.1,  we have to solve the following system of equations:
 \begin{equation}\left\{\begin{array}{l}
f_0=1,\mbox{ and } f_k=zf_{k-1}+zh_{k-1}, \quad k\geq 1,\\
g_k=z\sum\limits_{\ell\geq k+1} f_\ell+z\sum\limits_{\ell\geq k+1} h_\ell,\quad k\geq 0,\\
h_k=zf_k+zg_k+zh_k,\quad k\geq 0.\\
\end{array}\right.\label{equ2}
\end{equation}
Multiplying by $u^k$ the recursions in (\ref{equ2}) and summing over $k$, we have:
\begin{align*}
F(u)&=1+zuF(u)+zuH(u),\\
G(u)&=\frac{z}{u-1}(F(u)-F(1)+H(u)-H(1)),\\
H(u)&=zF(u)+zG(u)+zH(u).
\end{align*}

Notice that we have $(1-z)F(1)=1+zH(1)$ by considering the first equation. Now, setting $f1:=F(1)$ and solving these functional equations, we obtain
$$F(u)={\frac {{\it f1}\,u{z}^{2}-u{z}^{2}+zu+{z}^{2}-u-z+1}{{u}^{2}z+{z}^{2}
-u-z+1}},$$
$$G(u)=-{\frac {{\it f1}\,uz+{\it f1}\,z-zu-{\it f1}+1}{{u}^{2}z+{z}^{2}-u-z+
1}},\quad 
H(u)=-{\frac {z \left( {\it f1}\,uz-zu-{\it f1}+u+z \right) }{{u}^{2}z+{z}^
{2}-u-z+1}}.$$

In order to compute $f1$, we use the kernel method (see~\cite{ban, pro}) on $F(u)$. We can write the denominator (which is a polynomial in $u$ of degree 2), as $z(u-r)(u-s)$ with 
$$r={\frac {1+\sqrt {-4\,{z}^{3}+4\,{z}^{2}-4\,z+1}}{2z}}\quad\mbox{ and } \quad
s={\frac {1-\sqrt {-4\,{z}^{3}+4\,{z}^{2}-4\,z+1}}{2z}}.$$

Plugging $u=s$ (which has a Taylor expansion at $z=0$) in $F(u)z(u-r)(u-s)$, we obtain the equation 
${\it f1}\,s{z}^{2}-s{z}^{2}+zs+{z}^{2}-s-z+1 =0,$ which implies that $$f1=1+\frac{s-1}{z}.$$

Finally,  after simplifying by the factor $(u-s)$ in the numerators and denominators, we obtain 
 $$F(u)=\frac{r}{r-u},\quad G(u)=\frac{s-1}{z(r-u)},\quad\mbox{ and } \quad H(u)=\frac{s}{r-u},$$
which implies that
$$f_k=[u^k] F(u)=\frac{1}{r^k}, \quad g_k=[u^k] G(u)=\frac{s-1}{zr^{k+1}},\quad \mbox{ and }\quad
h_k=[u^k] H(u)=\frac{s}{r^{k+1}}.
$$

\begin{thm} The bivariate generating function for the total number of partial valley-less MAP  with respect to the  length  and the height of the end-point is given by
$$\mathit{Total}(z,u)=\frac{s}{z(r-u)},$$ and we have
$$[u^k] \mathit{Total}(z,u)=\frac{s}{zr^{k+1}}.$$
Finally, setting $t(n,k)=[z^n][u^k]\mathit{Total}(z,u)$, we have for $n\geq 2$, $k\geq 1$,  $$t(n,k)=t(n,k-1)+t(n-1,k)-t(n-1,k-2)-t(n-2,k),
$$
and setting $t_n:=t(n,0)$ and $t_{-1}=1$, then we have for $n\geq 2$, $$t_{n-1}=t_{n-2}+\sum\limits_{k=0}^{n-3}t_{k-1}t_{n-k-4}+\sum\limits_{k=2}^{n-1}\left(t_{k-1}-t_{k-2}\right)t_{n-k-2}.$$
\end{thm}
\noindent {\it Proof.} The first three equalities are immediately deduced from the previous results. For the last equality, the term $t_n$ satisfies the same recurrence relation as in  Theorem 1 (modulo shift of $n$) since the two sequences are equal modulo a shift.
\hfill $\Box$

\begin{cor} The generating function that counts the partial valley-less MAP with respect to the length is given by 
$$\mathit{Total}(z,1)=\frac{s}{z(r-1)}.$$
\end{cor}
The first few terms of the series expansion of $\mathit{Total}(z,1)$ are
$$1 + 2z + 5z^2 + 13z^3 + 34z^4 + 90z^5 + 242z^6 + 660z^7 + 1821z^8 + 5073z^9 + O(x^{10}),$$
which does not appear in \cite{oeis}.


\begin{cor} The generating function  that counts the partial valley-less MAP with respect to the length is given by 
$$\mathit{Total}(z,0)=\frac{s}{zr}.$$
\label{cor5}
\end{cor}
The first few terms of the series expansion of $\mathit{Total}(z,0)$ are  $$1 + z + 2z^2 + 5z^3 + 12z^4 + 29z^5 + 73z^6 + 190z^7 + 505z^8 + 1363z^9+ O(z^{10}),$$ which corresponds to a shift of the sequence \oeis{A152171} in \cite{oeis} counting Dyck paths of length $2n$ with no peaks of height $2$ (mod 3) and no valleys of height 1 (mod 3) (see also Section 3.1).

Let $\mathcal{T}$ be the infinite matrix $\mathcal{T}:=[t(n,k)]_{n, k\geq 0}$. The first few rows of the matrix $\mathcal{T}$ are
\[\mathcal{T}=
\left(
\begin{array}{ccccccccc}
 1 & 0 & 0 & 0 & 0 & 0 & 0 & 0 & 0 \\
 1 & 1 & 0 & 0 & 0 & 0 & 0 & 0 & 0 \\
 2 & 2 & 1 & 0 & 0 & 0 & 0 & 0 & 0 \\
 5 & 4 & 3 & 1 & 0 & 0 & 0 & 0 & 0 \\
 12 & 10 & 7 & 4 & 1 & 0 & 0 & 0 & 0 \\
 29 & 26 & 18 & 11 & 5 & 1 & 0 & 0 & 0 \\
 73 & 67 & 49 & 30 & 16 & 6 & 1 & 0 & 0 \\
 190 & 175 & 133 & 85 & 47 & 22 & 7 & 1 & 0 \\
 505 & 467 & 361 & 241 & 139 & 70 & 29 & 8 & 1 \\
\end{array}
\right).\]

\begin{figure}[H]
\centering
\includegraphics[scale=0.85]{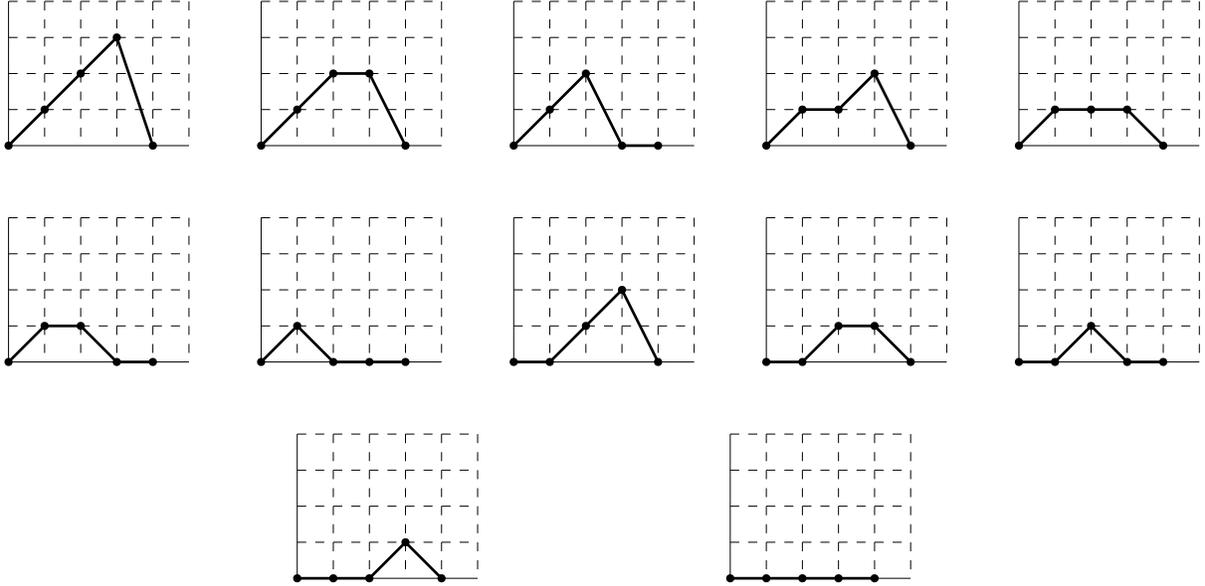}
\caption{The 12 $UU$-less PMAP of length $4$ ending at height $0$.} \label{fig3b}
\end{figure}

In Figure \ref{fig3b} we show the valley-less PMAP counted by $t(4,0)=12$. 

\begin{cor}
The matrix $\mathcal{T}=[t(n,k)]_{n, k\geq 0}$ is a Riordan array defined by 
$$\left((z^2-z+1)C\left(z(1-z+z^2)\right)^2,zC\left(z(1-z+z^2)\right)\right)$$
where $C(z)=\frac{1-\sqrt{1-4z}}{2z}$ is the generating function of the Catalan numbers
$c_n=\frac{1}{n+1}\binom{2n}{n}$.
\end{cor}
\noindent {\it Proof.}
 Indeed, we directly deduce the result from the following. 
\begin{align*}
[u^k]\mathit{Total}(z,u)&=\frac{s}{zr^{k+1}}=\frac{s}{zr}\cdot \frac{1}{r^k}=(z^2-z+1)C\left(z(1-z+z^2)\right)^2\cdot\left( zC\left(z(1-z+z^2)\right)\right)^k .
\end{align*}
\hfill $\Box$

From a similar argument as in Corollary \ref{corformula} we obtain the following result.
\begin{cor}
We have
$$t(n,k)=\sum_{j=0}^{n-k}\frac{k+2}{2(n-j)-k+2}\binom{2(n-j)-k+2}{n-k-j}a(n-k-j+1,j),$$
where $a(n,k)=(-1)^k\sum_{i=0}^n\binom ni \binom{n - i}{k - 2 i}$.
\end{cor}

\subsection{A bijective approach} 

Corollary~\ref{cor5} and Corollary~\ref{cor2} prove that the set of valley-less Motzkin paths with air pockets of length $n-1$ (ending on the $x$-axis) is equinumerous to the set of peak-less Motzkin paths with air pockets of length $n$, which is in one-to-one correspondence with the set $\mathcal{D}_n(2,1)$ of Dyck paths of length $2n$ with no peak at height 2 (mod 3) and no valley at height 1 (mod 3) (see a constructive bijection in Section~2.2). Below, we provide a bijection between valley-less MAP of length $n-1$ and peak-less MAP of length $n$. 

 Any  valley-less Motzkin path with air pockets is either of the form ($i$) $\epsilon$, ($ii$) $\alpha H$,  ($iii$) $U\alpha D$, ($iv$) $\beta H U\alpha D$, ($v$) $U\gamma D_k$, or ($vi$) $\beta H U\gamma D_k$, where $\alpha,\beta$ are valley-less MAP (possibly empty), and   $\gamma D_{k-1}$ is a valley-less MAP.  According to all these cases, we define the map $\phi$.
 
 \begin{defn} We  recursively define the map $\phi$ from  valley-less MAP of length $n-1$ to peak-less MAP of length $n$. Let $P$ be a valley-less MAP, we set:
$$
\phi(P)=\left\{\begin{array}{llr}
H&\text{if }P=\epsilon,&(i)\\
\phi(\alpha)H&\text{if }P=\alpha H,&(ii)\\
U\phi(\alpha)D&\text{if }P=U\alpha D,&(iii)\\
\phi(\beta)U\phi(\alpha)D&\text{if }P=\beta HU\alpha D,&(iv)\\
\phi(\alpha D_{k-1})^\sharp&\text{if }P=U\alpha D_k,&(v)\\
\phi(\beta)\phi(\alpha D_{k-1})^\sharp&\text{if }P=\beta HU\alpha D_k,&(vi)\\
\end{array}\right.
$$
\label{Def1b}
where the $\sharp$-operator maps a peak-less MAP of the form $\alpha D_{k-1}$ into the peak-less MAP $(\alpha D_{k-1})^\sharp=U\alpha D_k$.
\end{defn}
\medskip

Due to the recursive definition, the image  of valley-less MAP of 
length $n-1$ under $\phi$ is a peak-less MAP of length $n$. The recursive definition naturally induces that $\phi$ is a bijection. Using the bijection $\psi$ presented above, we can easily obtain a constructive bijection between valley-less MAP of length $n-1$ and Dyck paths of length $2n$ with no peak at height 2 (mod 3) and no valley at height 1 (mod 3).


\section{Partial $UU$-less Motzkin paths with air pockets}

In this section, we study partial Motzkin paths with air pockets of the first kind avoiding occurrences of $UU$. 
\subsection{Enumerative results}
In the same way as we done in the previous section,  we have to solve the following system of equations:
 
\begin{equation}\left\{\begin{array}{l}
f_0=1, f_1=z+zg_0+zh_0, \mbox{ and }
f_k=zg_{k-1}+zh_{k-1}, \quad k\geq 2,\\
g_k=z\sum\limits_{\ell\geq k+1} f_\ell+z\sum\limits_{\ell\geq k+1} h_\ell,\quad k\geq 0,\\
h_k=zf_k+zg_k+zh_k,\quad k\geq 0.\\
\end{array}\right.\label{equ1}
\end{equation}

Multiplying by $u^k$ the recursions in (\ref{equ1}) and summing over $k$, we have:
\begin{align*}
F(u)&=1+zu+zuG(u)+zuH(u),\\
G(u)&=\frac{z}{u-1}(F(u)-F(1)+H(u)-H(1)),\\
H(u)&=zF(u)+zG(u)+zH(u).
\end{align*}

Notice that we have $H(1)=(1+z)(F(1)-1)$ by considering the first and the last equation. Now, setting $f1:=F(1)$ and solving these functional equations, we obtain
$$F(u)=\frac{\mathit{f1} u \,z^{3}+2 \mathit{f1} u \,z^{2}+u^{2} z^{2}-u^{2} z -2 u \,z^{2}+2 u z +z^{2}-u -z +1}{u^{2} z^{2}+u \,z^{3}+u z +z^{2}-u -z +1}
,$$
$$G(u)=-\frac{z \left(\mathit{f1} u \,z^{3}+2 \mathit{f1} u \,z^{2}-u \,z^{3}+\mathit{f1} \,z^{2}-u \,z^{2}+\mathit{f1} z +u z -z^{2}-2 \mathit{f1} +2\right)}{u^{2} z^{2}+u \,z^{3}+u z +z^{2}-u -z +1}
,$$
$$H(u)=\frac{z \left(\mathit{f1} u \,z^{3}+2 \mathit{f1} u \,z^{2}-u \,z^{3}+\mathit{f1} \,z^{2}-u^{2} z -2 u \,z^{2}+2 \mathit{f1} z +u z -z^{2}-u -2 z +1\right)}{u^{2} z^{2}+u \,z^{3}+u z +z^{2}-u -z +1}
.$$

In order to compute $f1$,  we use the kernel method (see~\cite{ban, pro}) on $F(u)$. We can write the denominator (which is a polynomial in $u$ of degree 2), as $z^2(u-r)(u-s)$, with 
$$r={\frac {-{z}^{3}-z+1+\sqrt {{z}^{6}-2\,{z}^{4}+2\,{z}^{3}-3\,{z}^
{2}-2\,z+1}}{2{z}^{2}}} \
 \mbox{ and }$$
$$s={\frac {1-z-{z}^{3}-\sqrt {{z}^{6}-2\,{z}^{4}+2\,{z}^{3}-3\,{z}^{2}-
2\,z+1}}{2{z}^{2}}}
.$$

Plugging $u=s$ (which has a Taylor expansion at $z=0$) in $F(u)z^2(u-r)(u-s)$, we obtain the equation 
$\mathit{f1} s \,z^{3}+2 \mathit{f1} s \,z^{2}+s^{2} z^{2}-s^{2} z -2 s \,z^{2}+2 s z +z^{2}-u -z +1=0,$ which implies that $$f1=1+\frac{s-1}{z(z+2)}.$$
Finally,  after simplifying by the factor $(u-s)$ in the numerators and denominators, we obtain 
 $$F(u)=\frac{z-1}{z}+\frac{r}{z(r-u)},\quad G(u)=\frac{s+z}{r-u},\quad\mbox{ and } \quad H(u)=\frac{zr+1}{z(r-u)}-1,$$
which implies that $f_0=1$, $g_0=\frac{s+z}{r}$, $h_0=\frac{1}{zr}$ and
$$f_k=[u^k] F(u)=\frac{1}{zr^k}, \quad g_k=[u^k] G(u)=\frac{s+z}{r^{k+1}},\quad \mbox{ and }\quad
h_k=[u^k] H(u)=\frac{1+zr}{zr^{k+1}}.
$$


\begin{thm} The bivariate generating function for the total number of $UU$-less PMAP  with respect to the  length  and the height of the end-point is given by
$$\mathit{Total}(z,u)=\frac{1 + u z}{(r - u) z^2},$$ and we have $$[u^0] \mathit{Total}(z,u)=\frac{1}{z^2r},$$ and for $k\geq 1$
$$[u^k] \mathit{Total}(z,u)=\frac{rz+1}{z^2r^{k+1}}.$$
Finally, setting $t(n,k)=[z^n][u^k]\mathit{Total}(z,u)$, we have for $n\geq 2$, $k\geq 1$,
\begin{multline*}
    t(n,k)=t(n, k - 1) + t(n - 1, k) - t(n - 1, k - 1) - t(n - 2, k) \\- t(n - 2, k - 2) - t(n - 3, k - 1).
\end{multline*}
and setting $t_n:=t(n,0)$, then we have  $$t_n=t_{n-1}+t_{n-2}+t_{n-3}+\sum\limits_{k=2}^{n-2}t_{k-2}t_{n-k-2}+\sum\limits_{k=3}^{n-1}(t_{k-1}-t_{k-2})t_{n-k-1}.$$ 
\end{thm}
\noindent {\it Proof.} The first three equalities are immediately deduced from the previous results. Now, let us prove the last equality. 
 Any non-empty length $n$ $UU$-less MAP is of the form ($i$) $HP$ where $P$ is a $UU$-less MAP of length $n-1$, or ($ii$) $UDQ$ where $Q$ is a MAP of length $n-2$ avoiding $UU$, or ($iii$) $U Q D R$, where $Q,R$ are some MAP avoiding $UU$ such that the length of $Q$ lies into $[1,n-2]$ and $Q$ starts and ends with $H$, or ($iv$) $P Q$,  where $P=UP'D_{i}$, $i\geq 2$, and $P'D_{i-1}$ is a MAP  of length lying into $[3,n-1]$ and starting with $H$. The number of $P'D_{i-1}$ of a given length $k$ is the total number of $UU$-less MAP of length $k-1$ minus the total number of $UU$-less MAP of length $k-2$. Taking into account all these cases, we obtain the result.  

\hfill $\Box$

\begin{cor} The generating function that counts the partial $UU$-less MAP with respect to the length is given by 
$$\mathit{Total}(z,1)=\frac{1+z}{(r-1)z^2}.$$
\end{cor}

The first few terms of the series expansion of $\mathit{Total}(z,1)$ are
$$1 + 2z + 4z^2 + 9z^3 + 21z^4 + 50z^5 + 122z^6 + 302z^7 + 759z^8 + 1928x^{9}+O(x^{10}),$$
which does not appear in \cite{oeis}.

\begin{cor} The generating function that counts $UU$-less MAP with respect to the length is given by 
$$\mathit{Total}(z,0)=\frac{1}{z^2r}.$$
\label{cor11}
\end{cor}
The first few terms of the series expansion of $\mathit{Total}(z,0)$ are  $$1 + z + 2z^2 + 4z^3 + 9z^4 + 20z^5 + 47z^6 + 112z^7 + 274z^8 +  679x^9+O(z^{10}),$$ which corresponds to the sequence \oeis{A095980}  in \cite{oeis} counting Motzkin paths of length $n$ with no occurrences of $UHU$.

\begin{figure}[H]
\centering
\includegraphics[scale=0.85]{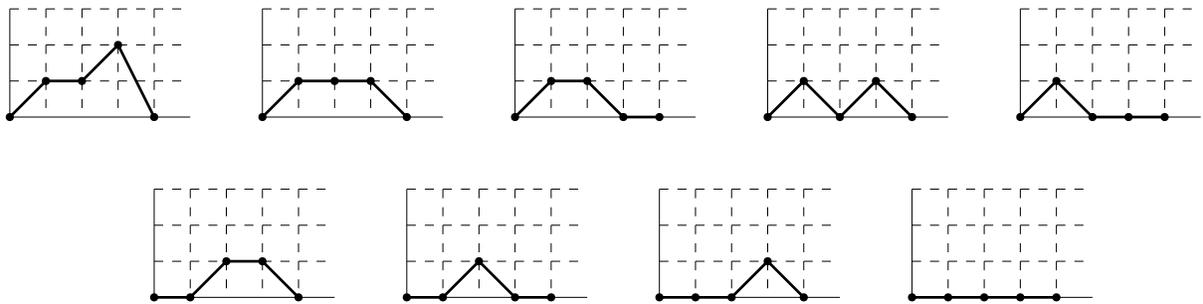}
\caption{The 9 $UU$-less PMAP of length $4$ ending at height $0$.} \label{fig4}
\end{figure}

Let $\mathcal{T}$ be the infinite matrix $\mathcal{T}:=[t(n,k)]_{n, k\geq 0}$. The first few rows of the matrix $\mathcal{T}$ are
\[\mathcal{T}=
\left(
\begin{array}{ccccccccc}
 1 & 0 & 0 & 0 & 0 & 0 & 0 & 0 & 0 \\
 1 & 1 & 0 & 0 & 0 & 0 & 0 & 0 & 0 \\
 2 & 2 & 0 & 0 & 0 & 0 & 0 & 0 & 0 \\
 4 & 4 & 1 & 0& 0 & 0 & 0 & 0 & 0 \\
 9 & 9 & 3 & 0 & 0 & 0 & 0 & 0 & 0 \\
 20 & 21 & 8 & 1 & 0 & 0 & 0 & 0 & 0 \\
 47 &50 & 21 & 4& 0 & 0 & 0 & 0 & 0 \\
 112& 121 & 55 & 13 & 1 & 0 & 0 & 0 & 0 \\
 274 & 298 & 143 & 39 & 5 & 0 & 0& 0& 0\\
\end{array}
\right).\]
In Figure \ref{fig4} we show the $UU$-less PMAP counted by $t(4,0)=9$. 

The matrix $\mathcal{G}$ is not a (proper) Riordan array. For this reason, we consider the matrix  $\mathcal{G}:=[g(n,k)]_{n\geq 0, k\geq 0}$, where
$$g(n,k)=\begin{cases}1, & \text{ if } n=k=0;\\
t(n+k-1,k), & \text{ if } n\geq 1.
\end{cases}$$

The first few rows of the matrix $\mathcal{G}$ are
\[\mathcal{G}=
\left(
\begin{array}{ccccccccc}
 1 & 0 & 0 & 0 & 0 & 0 & 0 & 0 & 0 \\
 1 & 1 & 0 & 0 & 0 & 0 & 0 & 0 & 0 \\
 1 & 2 & 1 & 0 & 0 & 0 & 0 & 0 & 0 \\
 2 & 4 & 3 & 1 & 0 & 0 & 0 & 0 & 0 \\
 4 & 9 & 8 & 4 & 1 & 0 & 0 & 0 & 0 \\
 9 & 21 & 21 & 13 & 5 & 1 & 0 & 0 & 0 \\
 20 & 50 & 55 & 39 & 19 & 6 & 1 & 0 & 0 \\
 47 & 121 & 143 & 113 & 64 & 26 & 7 & 1 & 0 \\
 112 & 298 & 372 & 319 & 203 & 97 & 34 & 8 & 1 \\
\end{array}
\right).\]

\begin{cor}
The matrix $\mathcal{G}=[g(n,k)]_{n, k\geq 0}$ is the  Riordan array defined by $(t(z),t(z)-1)$, where
$$t(z)=\frac{1 + z(1-z)^2 -\sqrt{(1 - 3 z + z^3) (1 + z + z^3)}}{2 z(1 - z + z^2)}.$$
\end{cor}
\noindent {\it Proof.}
It follows from the relation $t(z)=\frac{1}{zr}+1$.
\hfill $\Box$
\medskip

It seems difficult to obtain a close form for the coefficient $g(n,k)$ using the Lagrange Inversion Formula. Indeed,  $t(z)$ can be expressed in terms of $C(u)$, where $u=\frac{1}{4}(-z^6+2z^4-2z^3+3z^2+2z)$ is a polynomial of degree 6, which complicates the calculations.

\subsection{A bijective approach} 
Corollary~\ref{cor11} proves that the set of $UU$-less Motzkin paths with air pockets of length $n$ (ending on the $x$-axis) is equinumerous to the set of  Motzkin paths of length $n$ avoiding $UHU$. Below, we provide a bijection between these two sets. 

 Any  $UU$-less Motzkin path with air pockets is either of the form ($i$) $\epsilon$, ($ii$) $H\alpha$,  ($iii$) $UD\alpha $, ($iv$) $UHD\alpha$, ($v$) $UH\alpha H D\beta$, or ($vi$) $UH^k\gamma D_i \beta$, where $\alpha,\beta$ are $UU$-less MAP (possibly empty), and   $\gamma D_{i-1}$ is a $UU$-less MAP and $k\geq 1$.  According to all these cases, we define the map $\chi$.
 
 \begin{defn} We  recursively define the map $\chi$ from  $UU$-less MAP of length $n$ to $UHU$-less Motzkin paths of length $n$. Let $P$ be a $UU$-less MAP, we set:
$$
\chi(P)=\left\{\begin{array}{llr}
\epsilon&\text{if }P=\epsilon,&(i)\\
H\chi(\alpha)&\text{if }P=H\alpha,&(ii)\\
UD\chi(\alpha)&\text{if }P=UD\alpha,&(iii)\\
UHD\chi(\alpha)&\text{if }P=UHD\alpha,&(iv)\\
UHH\chi(\alpha) D\chi(\beta)&\text{if }P=UH\alpha H D\beta,&(v)\\
U\chi(\gamma D_{i-1})H^{k-1} D \chi(\beta)&\text{if }P=UH^k\gamma D_i \beta,&(vi)\\
\end{array}\right.
$$
\label{Def3}
\end{defn}
\medskip

\begin{figure}[H]
\centering
\includegraphics[scale=0.85]{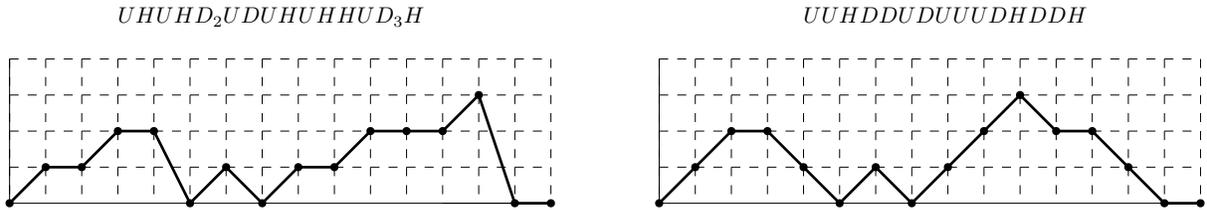}
\caption{$\chi(P)=UUHDDUDUUUDHDDH$.} \label{fig5}
\end{figure}

Due to the recursive definition, the image  of $UU$-less MAP of length $n$ by $\chi$ is a $UHU$-less Motzkin path  of length $n$. The recursive definition naturally induces that $\chi$ is a bijection. For instance, if  $P=UHUHD_2UDUHUHHUD_3H$ then we obtain $\chi(P)=\chi(UH^1UHD_2)\chi(UD)\chi(UHUHHUD_3)\chi(H)=UUHDDUDUUUDHDDH$ (see Figure \ref{fig5} for an illustration of this example).

\end{document}